\documentclass[12pt]{article}

\usepackage{t1enc,amssymb,amsbsy,amsthm,amsmath,amsfonts}

\newcommand{\RR}{\mathbb{R}}

\newcommand{\kk}{\mathcal{K}}
\newcommand{\ex}{\mathbf{E}}
\newcommand{\pr}{\mathbf{P}}
\newcommand{\prq}{\mathbf{Q}}
\newcommand{\rr}{\stackrel {d}{=}}
\newcommand{\indic}{\mathbf{1}}
\newcommand{\ps}{\widehat{\psi}_\alpha}
\newcommand{\fe}{\widehat{\varphi}_\alpha}

\newcommand{\feb}{\widehat{\varphi}_\beta}

\newtheorem{thm}{Theorem}[section]
\newtheorem{lem}[thm]{Lemma}
\newtheorem{prop}[thm]{Proposition}
\newtheorem{cor}[thm]{Corollary}
\newtheorem{exa}[thm]{Example}
\newtheorem{defi}[thm]{Definition}
\newtheorem{rema}[thm]{Remark}

\numberwithin{equation}{section}

\begin{document}

\title{On occupation times of stationary excursions}

\author{ \bf Marina Kozlova, \,
\bf Paavo Salminen
\vspace{0.4cm} \\
\it Department of Mathematics, \AA bo Akademi University, \\
\it \AA bo, FI-20500, Finland \\
e-mail addresses: mkozlova@abo.fi, phsalmin@abo.fi}

\date{\today}

\maketitle

\begin{abstract}
  In this paper excursions of a stationary diffusion in stationary state are studied.
  In particular, we compute
  the joint distribution of the occupation times $I^{(+)}_t$
  and $I^{(-)}_t$ above and below, respectively, the observed level at time $t$ during an
   excursion. We consider also
  the starting time $g_t$ and the ending time $d_t$
  of the excursion (straddling $t$) and discuss their relations to the Lévy measure
  of the inverse local time. It is seen that the pairs $(I^{(+)}_t, I^{(-)}_t)$ and
  $(t-g_t, d_t-t)$ are identically distributed. Moreover, conditionally on
  $I^{(+)}_t + I^{(-)}_t =v$, the variables $I^{(+)}_t$ and $I^{(-)}_t$ are uniformly distributed on $(0,v)$.
  Using the theory of the Palm measures, we derive an analoguous result for excursion bridges.
\end{abstract}



\section{Introduction} \label{intro}

In the literature there are many examples of cases such that
the occupation time for one diffusion is identical in law with the first hitting time
for another diffusion. For these see, e.g., \cite{ciesielskitaylor62,
biane85,  imhof86, legall85, yor92a,  donatimartinyor97, doneygrey89, pitmanyor03,
salminenyor04}.
The last one of these references contains a survey of known identities.

In this paper we prove a new kind of identity between hitting and occupation times for
stationary diffusions in stationary state. An example of such a diffusion is
reflecting Brownian motion on $\RR_+$ with drift $-\mu < 0$. As is well known, this process is
stationary having the exponential distribution with parameter $2\mu$ as its stationary
probability distribution. We take the whole of $\RR$ to be the
time axis and use, for a moment, the notation
$\{X_t\,:\, t\in\RR\}$ for this process. For fixed $t\in\RR$ let
\begin{equation} \label{defgd}
   g_t := \sup\{ s \leq t: X_s = 0\}, \quad d_t := \inf \{ s>t: X_s = 0\},
\end{equation}
and
$$
I^{(+)}_t:=\int_{g_t}^{d_t} \indic_{ \{X_s > X_t \} }\, ds,
\quad
I^{(-)}_t:=\int_{g_t}^{d_t} \indic_{ \{0<X_s \leq X_t \} }\, ds.
$$
Then from \cite{salminennorros01} it can be deduced
that
\begin{equation}
\label{result}
(I^{(+)}_t, I^{(-)}_t) \rr (t-g_t,d_t-t).
\end{equation}
For an explicit form of this distribution and its Laplace transform,
see Section 5, Example 5.1. The main result of the present  paper states
that (\ref{result}) is valid for {\sl all} stationary diffusions in stationary state.
By time reversal, the variables $I^{(+)}_t, I^{(-)}_t, t-g_t,$ and $d_t-t$ are
identically distributed.
Our proof of (\ref{result}) is purely computational and does not, unfortunately,
provide
any probabilistic explanation of the identity. In the case of the reflecting Brownian
motion with drift we have an alternative, more probabilistic, proof based on the Ray--Knight
theorems. It is fairly easy to deduce also in the general case using the symmetry and time
reversal as is done in \cite{salminennorros01}, that the expectations
of $I^{(+)}_t, I^{(-)}_t, t-g_t,$ and $d_t-t$ are equal.

It is of interest to recall the identity due to \cite{biane85} and 
\cite{imhof86} because the variable $d_t-t$ can also be found in this identity.
Let $\{B^{(\mu)}_t\, :\, t\geq 0\}$ be a
Brownian motion with drift $\mu>0$ starting from 0 and
$$
H_x(B^{(\mu)}):=\inf\{ t\,:\, B^{(\mu)}_t=x\}
$$
the first hitting time of the level $x>0$. Then the Biane--Imhof identity states that
\begin{equation}
\label{i}
\int_0^\infty {\bf 1}_{\{B^{(\mu)}_s< 0\}}\, ds
\quad{\mathop=^{\rm{(d)}}}\quad
H_\lambda(B^{(\mu)}),
\end{equation}
where $\lambda$ is exponentially (with parameter $2 \mu$) distributed random
variable independent of $B^{(\mu)}$.
 By spatial homogeneity,
we have
$$
H_\lambda(B^{(\mu)})
\quad{\mathop=^{\rm{(d)}}}\quad
d_t-t
$$
since the distribution of $X_t$ is exponential with parameter $2\mu$.

From the identity (\ref{result}) using the theory of Palm measures applied for excursions
we derive a new interesting result for excursion bridges. To explain this,
let $X^{(0,l,0)}$ be an excursion bridge of length $l$ for excursions from 0 to 0 for an
arbitrary diffusion $X.$ Further, let $U$ be
a random variable having the uniform distribution
on $(0,l)$ and assume that $U$ is independent of $X^{(0,l,0)}.$ Then the occupation times
$$
I^{(l,+)}:=\int_{0}^{l} \indic_{ \{X^{(0,l,0)}_s > X^{(0,l,0)}_U \} }\, ds,
\quad
I^{(l,-)}:=\int_{0}^{l} \indic_{ \{X^{(0,l,0)}_s < X^{(0,l,0)}_U \} }\, ds
$$
are identically uniformly distributed on $(0,l)$. For a
standard Brownian excursion, i.e., for the 3-dimensional Bessel bridge,
this fact can also be explained via Vervaat's
 path transformation.

The paper is organized so that there are two main sections
and two shorter section with remarks and examples. The first main section is devoted
to the identity (\ref{result}) and is divided into four subsections. In the first subsection
some preliminaries on stationary diffusions are given. After this the joint
distribution of $g_0$ and $d_0$ is derived. In particular, it is seen that
this distribution belongs to a special class of two-dimensional distributions
characterized by the fact that the density (when it exists)
is  a function of the sum of the arguments only. The proof of the main identity
(\ref{result}) is presented in Section $\ref{identity}$. After this, in Section
$\ref{subrbm}$ we give an alternative
proof of (\ref{result}) for RBM with drift based on the Ray--Knight theorems.
In the second
main section the corresponding identity for the excursion bridges is discussed.
Here we start with by relating the distribution of $(-g_0,d_0)$ to the L\'evy measure
of the inverse local time preparing in this way the connection of
the It\^{o} excursion measure and the Palm measure. In Section $\ref{subspectral}$
we make some observations concerning the spectral representations of $d_0$ and
$d_0-g_0$. The result concerning excursion bridges is proved in Section $\ref{palm}$ by
utilizing the connections of the It\^{o} excursion measure and the Palm measure.
In Section $\ref{nullrec}$ we discuss the null recurrent case and, finally,
in Section $\ref{subexamples}$ we give some examples.

\section{Stationary excursions straddling $t$} \label{results}

\subsection{Preliminaries} \label{prel}
Let $X^\circ = \{X^\circ_t: t\geq 0\}$ be a one-dimensional recurrent conservative
diffusion living in the interval $I$. We assume that
$I=[0,b]$, $b < + \infty$, or $I=[0, b)$, $b \leq +\infty$, where $0$ is a reflecting
boundary and $b$ is either reflecting or natural or entrance-not-exit.
 For background on
one-dimensional diffusions we refer to  \cite{itomckean74} and
\cite{borodinsalminen02}.

The diffusion $X^\circ$ is characterized in $I$ by its
scale function $S(x)$ and speed measure $m(dx)$.
Suppose, moreover, that $X^\circ$ is positively recurrent, that is $M:= m\{I\} < \infty$.
Let $p(t;x,y)$ be the symmetric transition
density of $X^\circ$ with respect to $m(dx)$ and let
$$\mathcal{G} = \frac{d}{dm} \frac{d}{d S}$$
be the infinitesimal generator of $X$. Assume that
$m(dx) = m(x) dx$, $S(dx)= S'(x) dx$ such that $m(x)$ and $S'(x)$ are continuous and
positive.
The Green function is defined as
$$G_\alpha (x,y):= \int_0^\infty e^{ -\alpha t} p(t;x,y) \, dt.$$
Let $\pr_x$ and $\ex_x$ denote the probability measure and the expectation, respectively,
associated with
$X^\circ$ started at $x$ and $H_y(X^\circ)$ denote the first hitting time of $y$ for $X^\circ$.
Recall from  \cite{itomckean74} p. 129
 that
\begin{equation}\label{hit}
 \ex_x \left( e^{-\alpha H_y(X^\circ)}\right) = \frac{G_\alpha(x,y)}{G_\alpha(y,y)}.
\end{equation}

Let now $\{X_t^{(1)} : t \geq 0\}$ and $\{X_t^{(2)}: t \geq 0\}$ be two copies of $X^\circ$
such that $X_0^{(1)}=X_0^{(2)}$ with the common law
$\widetilde m(dx):= m(dx) /M$ but let $X^{(1)}$ and $X^{(2)}$ be otherwise
independent. Define for $t \in \RR$,
$$X_t:=\left \{
    \begin{array}{cc}
       X_t^{(1)}, & t \geq 0, \\
       X_{-t}^{(2)}, & t \leq 0.
     \end{array}
     \right.$$
The process $X=\{X_t: t \in \RR\}$ is called a stationary diffusion in stationary state
living in $I$ and having the generator $\mathcal{G}$. Notice that the law of
$X_t$ is $\widetilde m$ for every $t \in \RR$.


\subsection{Joint distribution of $g_0$ and $d_0$} \label{sub1}

For the process $X$ defined above let $g_0$ and $d_0$ be as  in $(\ref{defgd})$.
In this section we derive the joint distribution of $-g_0$ and $d_0$ and show that
the pair $(-g_0, d_0)$
is an element of a special class $\kk$ of
two-dimensional random variables studied in 
\cite{salminenvallois03} and \cite{kozlovasalminen03}. We start with
the definition and some properties of this class $\kk$.

\begin{defi}\label{defk}
\rm A two-dimensional random variable $(\xi_1,\xi_2)$, where both $\xi_1$ and
$\xi_2$ are non-negative, belongs to $\kk$ if
$$(\xi_1,\xi_2) \rr(U, V-U),$$
where $V$ is an arbitrary non-negative random variable and the conditional
distribution of $U$ given $V$
is uniform on $(0,V)$.
\end{defi}

\begin{prop} \label{propk} \bf (Properties of the elements in $\kk$) \it \\
 1) \; Let $(\xi_1,\xi_2) \in \kk$. Then $\xi_1 \rr \xi_2$, $\xi_1+\xi_2 \rr V$, and the
 density of $\xi_1$ (and $\xi_2$) exists and is given
by $$F'_{\xi_1}(x) = \int_{[x, \infty)} v^{-1} F_V(dv).$$
If the density $F'_V(v)=:p(v)$ exists then
\begin{equation} \label{eqk}
\pr (\xi_1 \in dx, \xi_2 \in dy) = \frac{p(x+y)}{x+y}\,dx\,dy.
\end{equation}
2) \; Let $\xi_1$ and $\xi_2$ be non-negative random variables. Then $(\xi_1,\xi_2) \in \kk$
 if and only if
 for all $\alpha \ne \beta$,
   \begin{equation} \label{formulakl}
   \ex \left(e^{-\alpha \xi_1-\beta \xi_2} \right) = \frac{1}{\alpha-\beta} \int_\beta^\alpha
   \ex \big(e^{-\gamma (\xi_1+\xi_2)} \big) \,d \gamma.
   \end{equation}
\end{prop}
\noindent
\textit{Proof:} For the proof, see \cite{salminenvallois03} and \cite{kozlovasalminen03}. \qed

\begin{prop} \label{proposition1}
  The pair $(-g_0, d_0)$ belongs to the class $\kk$
  with the joint Laplace transform given by
  \begin{equation} \label{jointlap}
    \ex \left( e^{-\alpha d_0 + \beta g_0} \right) =
    \frac{1}{M(\alpha-\beta)} \left( \frac{1}{G_{\alpha} (0,0)}
    -\frac{1}{G_{\beta}(0,0)} \right).
  \end{equation}
  Moreover, the joint density of $(-g_0,d_0)$ exists and is
  given by
\begin{eqnarray} \label{joint}
  &&\pr(d_0 \in dt,  -g_0 \in ds) \nonumber \\
   && \hspace{1cm} = \frac{1}{M}\left ( \frac{d}{dS(y_1)} \frac{d}{dS(y_2)}
    \widehat p(t + s; y_1, y_2) \right) \Big|_{y_1, y_2=0+}\, dt \, ds,
\end{eqnarray}
where $\widehat p(t;x,y)$ is the transition density of the  process
$\widehat X$ obtained from $X^\circ$ by killing  $X^\circ$ when it hits zero.
\end{prop}
\noindent
\textit{Proof:}
Using conditional independence of $X^{(1)}$ and $X^{(2)}$, formula $(\ref{hit})$,
and the Chapman-Kolmogorov equation, we write
\begin{eqnarray*}
& & \ex \left( e^{-\alpha d_0 +\beta g_0} \right)
 = \int_I \widetilde m(dx) \, \ex_x \big(e^{-\alpha H_0(X^{(1)})} \big)
  \ex_x \big(e^{-\beta H_0(X^{(2)})}
  \big) \\
  & & \hspace{0.4cm} =
      \int_I  \widetilde m(dx) \, G_{\alpha} (x,0) G_{\beta} (x,0) \,
      /(G_{\alpha} (0,0) G_{\beta} (0,0))\\
  & & \hspace{0.4cm} =
      \int_I \widetilde m(dx) \int_0^\infty dt \, e^{-\alpha t} p(t;x,0) \\
   & & \hspace{1.7cm} \times
       \int_0^\infty ds \, e^{-\beta s} p(s;x,0)
       /(G_{\alpha} (0,0) G_{\beta} (0,0)) \\
& & \hspace{0.4cm} =
      \int_0^\infty dt\int_0^\infty ds \, e^{-\alpha t- \beta s} \\
  & & \hspace{1.7cm} \times
      \int_I  \widetilde m(dx)\, p(t;x,0) p(s;x,0)  /(G_{\alpha} (0,0) G_{\beta}(0,0)) \\
 & & \hspace{0.4cm} =
      \frac{1}{M}\int_0^\infty dt \int_0^\infty ds \,e^{-\alpha t- \beta s}
       p(t+s;0,0)
   \,   /(G_{\alpha} (0,0) G_{\beta}(0,0)) \\
 & & \hspace{0.4cm} = \frac{1}{M(\alpha-\beta)} \left( \frac{1}{G_{\alpha} (0,0)}
   -\frac{1}{G_{\beta}(0,0)} \right).
\end{eqnarray*}
To show that $(-g_0,d_0) \in \kk$ we use Proposition $\ref{propk}$. Letting $\alpha \to \beta$
in $(\ref{jointlap})$, it is seen that the derivative
$\frac{d}{d \alpha}\big(\frac{1}{G_\alpha(0,0)}\big)$ exists and
\begin{eqnarray*}
\ex \left( e^{-\alpha d_0 +\beta g_0} \right)
   & = &\frac{1}{M(\alpha-\beta)} \int_\beta^\alpha
   \frac{d}{d \gamma}\Big( \frac{1}{G_\gamma(0,0)}\Big) \, d \gamma \\
   & = & \frac{1}{M(\alpha-\beta)} \int_\beta^\alpha
   \ex \left( e^{-\gamma (d_0 -g_0)}\right) d \gamma.
\end{eqnarray*}
To compute the joint density recall
the formula (see  \cite{itomckean74}, p. 154)
\begin{equation} \label{itoh}
  \pr_x(H_0(X^\circ) \in dt) /dt = \frac{d \widehat p(t;x,y)}{d S(y)}\Big|_{y=0+}
  =:\widehat p^+(t;x,0)
\end{equation}
and, hence,
\begin{eqnarray*}
   &&  \pr(d_0 \in dt,  -g_0 \in ds) / dt\, ds \\
   & & \hspace{0.7cm} = \int_I \widetilde m(dx)
   \widehat p^+(t;x,0) \widehat p^+(s;x,0)\\
    & & \hspace{0.7cm} = \frac{1}{M}\left ( \frac{d}{d S(y_1)}
     \frac{d}{dS(y_2)}
    \int_I m(dx)\, \widehat p(t; x,y_1) \widehat p(s; x,y_2)
    \right) \Big|_{y_1, y_2=0+} \\
    & & \hspace{0.7cm} = \frac{1}{M}\left ( \frac{d}
    {d S(y_1)} \frac{d}{d S(y_2)}
    \widehat p(t + s; y_1, y_2) \right) \Big|_{y_1, y_2=0+},
\end{eqnarray*}
where it is used
that $\widehat p(t;x,y)$ is continuous  in $t$, $x$, $y$ and that $\widehat p^+(t;x,y)$
is continuous in $t$, $x$ and right continuous in $y$ (see
 \cite{itomckean74} p. 149). \qed

\begin{rema} \label{ghat}
\rm
To find the density $(\ref{joint})$ we used formula $(\ref{itoh})$ and the
Chapman-Kolmogorov equation.
There is also an alternative approach. For this,
let $\widehat G_\alpha (x,y)$ be the Green function for the killed process
$\widehat X$. Then (cf.  \cite{csakifoldes87} p. 185)
 \begin{equation} \label{eqg}
   \frac{d}{dS(x)} \frac{d}{d S(y)} \widehat G_\alpha (x,y) \Big{|}_{x,y=0+}
   = - \frac{1}{G_\alpha(0,0)}
 \end{equation}
and we have
\begin{eqnarray*}
  & & \frac{1}{\alpha-\beta} \left( \frac{1}{G_\alpha(0,0)} -
  \frac{1}{G_\beta(0,0)}\right) \\
  & & \hspace{1.2cm} = \frac{1}{\alpha-\beta} \frac{d}{dS(x)} \frac{d}{dS(y)}
  \left( \widehat G_\beta(x,y) - \widehat G_\alpha(x,y)\right) \Big{|}_{x,y=0+} \\
  & & \hspace{1.2cm} = \frac{1}{\alpha-\beta} \int_0^\infty dt\, \left(
  e^{-\beta t} -e^{-\alpha t}\right) \frac{d}{dS(x)} \frac{d}{dS(y)}
  \widehat p(t; x,y) \big{|}_{x,y=0+}\\
  & & \hspace{1.2cm} = \int_0^\infty \int_0^\infty e^{-\alpha t -\beta s}
  \frac{d}{dS(x)} \frac{d}{dS(y)} \widehat p(t+s; x,y) \big{|}_{x,y=0+} \, dt \, ds.
\end{eqnarray*}
\end{rema}

\subsection{Main identity in law} \label{identity}

The first main result of the
paper is presented in the following theorem.

\begin{thm} \label{theorem1}
 Let $X$ be a stationary diffusion in stationary state. Then for all
 $t \in \RR$, the two-dimensional random variables $(t-g_t, d_t-t)$ and
 $(I^{(+)}_t, I^{(-)}_t)$  are identically
 distributed and, from Proposition $\ref{proposition1}$,
  \begin{equation} \label{jointlapii}
    \ex \big( e^{-\alpha I^{(+)}_t - \beta I^{(-)}_t} \big) =
    \frac{1}{M(\alpha-\beta)} \left( \frac{1}{G_{\alpha} (0,0)}
    -\frac{1}{G_{\beta}(0,0)} \right),
  \end{equation}
 \begin{eqnarray} \label{jointii}
  && \pr \big(I^{(+)}_t \in dr,  I^{(-)}_t \in ds \big)  \nonumber\\
    && \hspace{1cm} = \frac{1}{M}\left ( \frac{d}{dS(y_1)} \frac{d}{dS(y_2)}
    \widehat p(s + r;\, y_1, y_2) \right) \Big|_{y_1, y_2=0+}\, dr \, ds.
 \end{eqnarray}
\end{thm}
\noindent
\textit{Proof:}
We prove the theorem only for the case when the state space is $I=[0, +\infty)$ and leave
other cases to the
reader.

By the stationarity, it suffices to consider the case $t=0$.
To simplify the notation, let
$I^{(+)}:=I^{(+)}_0$ and $I^{(-)}:=I^{(-)}_0$.
For $y \in (0, \infty)$
introduce
$$u(x):=\ex_x \left( \exp\Big( -\alpha \int_0^{H_0(X)}
\indic_{ \{ 0 \leq X_s \leq y \} } ds  - \beta \int_0^{H_0(X)} \indic_{ \{ X_s >y \}} ds
\Big)\right).$$
From the Feynman-Kac formula it follows that
$u(x)$ is the unique bounded smooth solution
of the equation
$$\mathcal{G} u(x) = \left\{
   \begin{array}{ll}
    \alpha \,u(x), & 0 <x <y, \\
    \beta \,u(x), & x >y
   \end{array}
\right.$$
with the boundary condition $u(0)=1$.
The function  $u(x)$, therefore, has the following form
$$u(x):=\left \{
    \begin{array}{ll}
       \vspace{0.2cm}
       A \ps(x) +  \fe(x) /\fe(0), & x \leq y, \\
      B \feb(x), & x \geq y,
     \end{array}
     \right.$$
for some constants $A$ and $B$,
where
$\widehat \psi_\alpha$ and $\widehat \varphi_\alpha$ are the
increasing and decreasing fundamental solution, respectively, of the
generalized differential equation
\begin{equation} \label{ur1}
  \mathcal{G} u = \alpha u
\end{equation}
such that $\fe$ is bounded and the "killing" condition holds
\begin{equation} \label{b2}
   \ps(0) =0.
\end{equation}
Both $u$ and its derivative $u^+$ should be continuous at $y$, i.e.,
\begin{equation} \label{AB}
\left\{
   \begin{array}{cc}
    \vspace{0.2cm}
    A \ps(y) +  \fe(y) /\fe(0) = B \feb(y) \\
    A \ps^+(y) + \fe^+(y) / \fe(0) = B \feb^+(y).
    \end{array}
    \right.
\end{equation}
From $(\ref{AB})$,
\begin{eqnarray*}
B & = & \frac{\ps^+(y) \fe(y) - \ps(y) \fe^+(y)}{\fe(0) \big(\ps^+(y) \feb(y) -
\ps(y) \feb^+(y) \big)}\\
& = & \frac{\ps^+(0)}{\ps^+(y) \feb(y) - \ps(y) \feb^+(y)},
\end{eqnarray*}
where it is used that
the Wronskian (a constant)
\begin{equation} \label{wr}
  \widehat w_\alpha:= \ps^+(x) \fe(x) - \ps(x) \fe^+(x)
\end{equation}
equals  $\ps^+(0) \fe(0)$.
Consequently, at point $y$,
$$u(y) = \frac{\ps^+(0) \feb(y)}{\ps^+(y) \feb(y) - \ps(y) \feb^+(y)}.$$
The joint Laplace transform of $I^{(-)}$ and $I^{(+)}$ is given then as
\begin{eqnarray} \label{i-i+}
  &&\ex \big( e^{-\alpha I^{(-)} -\beta I^{(+)}} \big) =  \int_0^\infty \widetilde m(dy)\, (u(y))^2 \nonumber\\
   && \hspace{0.6cm }= \big(\ps^+(0) \big)^2 \int_0^\infty \widetilde m(dy)
   \left(\frac{\feb(y)}{\ps^+(y) \feb(y) - \ps(y) \feb^+(y)}
   \right)^2.
\end{eqnarray}
Let
$$r(y):= \ps^+(y) \feb(y) - \ps(y) \feb^+(y)$$
denote the denominator in the
expression for $u(y)$
and notice that
\begin{eqnarray} \label{a-b}
  \frac{d}{dm} r(y) & = & \frac{d}{dm} \left( \feb(y) \frac{d}{dS} \ps(y)
  - \ps(y) \frac{d}{dS} \feb(y)\right) \nonumber\\
 & = & (\alpha -\beta) \ps(y) \feb(y),
\end{eqnarray}
where we used that $m(dx)= m(x) dx$ and $S(dx) = S'(x) dx$
and that $\ps(x)$ and $\feb(x)$ satisfy the equations
$$\mathcal{G} \ps = \alpha \ps \quad \text{and} \quad
\mathcal{G} \feb = \beta \feb,$$
respectively.
The equality $(\ref{a-b})$ is the key to obtain
\begin{eqnarray*}
  \ex \big( e^{-\alpha I^{(-)} -\beta I^{(+)}} \big) & = & \big(\ps^+(0) \big)^2 \int_0^\infty \widetilde m(dy)
  \frac{\big(\feb(y)\big)^2}{\big(r(y)\big)^2} \\
   & = & \frac{\big(\ps^+(0)\big)^2}{M(\alpha-\beta)} \int_0^\infty m(dy)
   \frac{\feb(y)}{\ps(y)} \frac{\frac{d}{dm} r(y)}{\big(r(y)\big)^2} \\
   & = & \frac{\big(\ps^+(0)\big)^2}{M(\alpha-\beta)}  \lim_{\varepsilon \to 0} \left\{
    \int_\varepsilon^\infty m(dy)
   \frac{\feb(y)}{\ps(y)} \frac{\frac{d}{dm} r(y)}{\big(r(y)\big)^2} \right\}.
\end{eqnarray*}
Let $F(\varepsilon)$ denote the integral above.
Using that
$\feb(x)$ is bounded and $\ps(x) \to +\infty$ as $x \to +\infty$,
we have
\begin{eqnarray*} F(\varepsilon)  & = & \int_\varepsilon^\infty m(dy)
   \frac{\feb(y)}{\ps(y)} \frac{d}{dm} \left( -\frac{1}{r(y)}\right) \\
    & = & \int_\varepsilon^\infty m(dy)
   \frac{d}{dm} \left( \frac{1}{r(y)} \right) \int_y^\infty S(dx) \frac{d}{dS}
   \left( \frac{\feb(x)}{\ps(x)}\right) \\
& = &
\int_\varepsilon^\infty S(dx) \frac{d}{dS} \left( \frac{\feb(x)}{\ps(x)}\right)
\int_\varepsilon^x m(dy) \frac{d}{dm} \left( \frac{1}{r(y)}\right) \\
& = &  \int_\varepsilon^\infty S(dx) \frac{d}{dS} \left(
\frac{\feb(x)}{\ps(x)}\right) \left( \frac{1}{r(x)}
 -\frac{1}{r(\varepsilon)}\right).
\end{eqnarray*}
Since
$$r(x) = - \big(\ps(x) \big)^2 \frac{d}{dS} \left( \frac{\feb(x)}{\ps(x)}\right),$$
we obtain
$$F(\varepsilon)  = - \int_\varepsilon^\infty
\frac{S(dx)}{\big(\ps(x)\big)^2} + \frac{\feb(\varepsilon)}{\ps(\varepsilon)
r(\varepsilon)}. $$
By definition ($\ref{wr}$) of the Wronskian,
$$\frac{1}{\widehat{w}_\alpha} \frac{d}{dS} \left( \frac{\fe(x)}{\ps(x)} \right)=
-\frac{1}{\big(\ps(x) \big)^2},$$
and, hence,
\begin{eqnarray*}
F(\varepsilon) & = & -\frac{1}{\widehat w_\alpha}
\frac{\fe(\varepsilon)}{\ps(\varepsilon)} + \frac{\feb(\varepsilon)}
{\ps(\varepsilon) \big(\ps^+(\varepsilon) \feb(\varepsilon) -
\ps(\varepsilon) \feb^+(\varepsilon) \big)}  \\
& = & \frac{1}{\ps(\varepsilon)}
\left(\frac{-\fe(\varepsilon)}{\ps^+(\varepsilon) \fe(\varepsilon)
- \ps(\varepsilon) \fe^+(\varepsilon)}
+ \frac{\feb(\varepsilon)}
{ \ps^+(\varepsilon) \feb(\varepsilon) -
\ps(\varepsilon) \feb^+(\varepsilon)} \right) \\
& = & \frac{\fe(\varepsilon) \feb^+(\varepsilon) -\feb(\varepsilon)
\fe^+(\varepsilon)}{\widehat w_\alpha \big(\ps^+(\varepsilon) \feb(\varepsilon) -
\ps(\varepsilon) \feb^+(\varepsilon) \big)}.
\end{eqnarray*}
Consequently, using boundary condition $(\ref{b2})$,
\begin{eqnarray} \label{jointlapi}
  \ex \big( e^{-\alpha I^{(-)} -\beta I^{(+)}} \big) & = & \frac{\big(\ps^+(0)\big)^2}{M(\alpha-\beta)}
  \lim_{\varepsilon \to 0} F(\varepsilon) \nonumber \\
  & = & \frac{\big(\ps^+(0)\big)^2}{M(\alpha-\beta)}\, \frac{\big(\fe(0) \feb^+(0) -
  \feb(0) \fe^+(0) \big)}{\widehat w_\alpha \ps^+(0) \feb(0)} \nonumber \\
  & = &\frac{1}{M(\alpha - \beta)}
  \left(\frac{\feb^+(0)}{\feb(0)} - \frac{\fe^+(0)}{\fe(0)} \right).
\end{eqnarray}
Finally,
the Green function of $X$ can be represented as
$$G_\alpha(x,y) = w_\alpha^{-1} \psi_\alpha(x) \varphi_\alpha(y), \quad x \leq y$$
 (see, e.g.,
\cite{borodinsalminen02} p. 19),
where $\psi_\alpha$,  $\varphi_\alpha$ are the fundamental solutions of
$(\ref{ur1})$ such that $\varphi_\alpha$ is bounded and $\psi_\alpha^+(0)=0$ (condition for
the reflection).
Since $\fe = \varphi_\alpha$, we have
$$\frac{1}{G_\alpha(0,0)} = - \frac{\varphi^+_\alpha(0)}{\varphi_\alpha(0)}=
-\frac{\fe^+(0)}{\fe(0)}.$$
Thus, the right-hand sides of $(\ref{jointlapi})$ and
$(\ref{jointlap})$ are equal and the proof is complete.
\qed


Although we do not have any probabilistic explanation for
the equality in law in Theorem $\ref{theorem1}$ (but for $RBM$ with drift, see the
next section), the following lemma (cf. \cite{salminennorros01}, Proposition 3.11 p. 330)
explains why the variables $I^{(+)}_t$, $I^{(-)}_t$,
$t-g_t$, and $d_t-t$ all have the same expectation.
\begin{lem}\label{prop311}
  For a given $s>0$ the events $\{X_s > X_0\}$ and $\{d_0>s\}$ are independent.
\end{lem}
\textit{Proof:} \rm
By the symmetry of the transition densities,
\begin{eqnarray} \label{form311}
  & & \pr(X_s > X_0 \,;\, s < d_0) = \int_0^b \frac{m(dx)}{M} \pr_x (X_s^{(1)} >x\,;\, s
  < H_0(X^{(1)}))
  \nonumber \\
  & & \hspace{2cm} = \frac{1}{M} \int_0^b m(dx) \int_x^b m(dy) \, \widehat p(s; x,y)  \nonumber\\
  & & \hspace{2cm} = \int_0^b \frac{m(dy)}{M} \int_0^y m(dx) \,\widehat p(s;y,x) \nonumber \\
  & & \hspace{2cm} = \pr(X_s < X_0 \,;\, s< d_0).
\end{eqnarray}
Noting that, by the reversibility and the stationarity,
$$\pr(X_s>X_0) = \pr(X_{-s} > X_0) = \pr(X_0>X_s) = 1/2$$
 completes the proof. \qed

\begin{cor} \label{corex}
  The expectations of $I^{(+)}$,  $I^{(-)}$, and $d_0$ are equal.
\end{cor}
\noindent
\textit{Proof:}
By $(\ref{form311})$,
\begin{eqnarray*}
  & & \ex \left( I^{(+)} \right) = 2\, \ex \Big( \int_0^{d_0} \indic_{ \{X_s > X_0 \}}\, ds \Big)
  = 2 \int_0^\infty \pr(X_s> X_0 \,;\, s < d_0) \, ds \\
  & & \hspace{0.7cm} = \int_0^\infty \pr(d_0 > s) \, ds = \ex (d_0).
  \end{eqnarray*}
Similarly, we can show that $\ex \left( I^{(-)} \right) = \ex (d_0)$.
\qed

\subsection{Probabilistic explanation of the identity when $X$ is a RBM with drift} \label{subrbm}

In the case when  $X$ is a stationary reflecting Brownian motion with drift $-\mu <0$ we
prove the fact that $I^{(+)} \rr I^{(-)} \rr d_0$  by considering the local
time process of an excursion instead of applying the Feynman-Kac
formula as in the proof of Theorem $\ref{theorem1}$.

Let $\{L(t,y): t \geq 0, y \geq 0\}$ denote the local time of $\{X^\circ_s : s \geq 0\}$ up to $t$ at the
level $y$ with respect to the Lebesgue measure, i.e., we have
\begin{equation} \label{occup}
 \int_0^t g(X^\circ_s) \, ds = \int_E g(x) L(t,x) \,dx
\end{equation}
for any non-negative Borel function $g$.
The following Ray--Knight theorem (see 
\cite{borodinsalminen02} pp. 90--91) describes the behaviour of the local
time process $L$ up to $H_0(X)$. Let $Z^{(n, 2 \mu)}$ denote the squared radial
part of an $n$-dimensional Ornstein--Uhlenbeck
process with parameter $\mu$ and recall that the generator of $Z^{(n, 2 \mu)}$ is
$$2 z \frac{d^2}{dz^2} + (n - 2 \mu z) \frac{d}{dz}.$$

\begin{thm} \label{theorem4}
  Conditionally on $X_0=x$,
$$\{ L(H_0(X), y): 0 \leq y \leq x\} \rr \{Z_y^{(2, 2 \mu)}: 0 \leq y \leq x\},$$
$$\{ L(H_0(X), x+y):y \geq 0\} \rr \{Z_y^{(0, 2 \mu)}: y \geq 0\},$$
where the process $Z^{(0, 2 \mu)}$ is started from the position of
$Z^{(2, 2 \mu)}$  at time $x$ but otherwise  $Z^{(2, 2 \mu)}$ and $Z^{(0, 2 \mu)}$
are independent.
\end{thm}

For the stationary excursion
$\{X_t: g_0 < t < d_0\}$ straddling zero
introduce its  total local time process by
$$L^{(e)}(y):= L^{(1)} (H_0(X^{(1)}),y)+ L^{(2)}(H_0(X^{(2)}),y), \quad y \geq 0,$$
where $L^{(1)} (H_0(X^{(1)}),y)$ and $L^{(2)} (H_0(X^{(2)}),y)$ are the local times at $y$
up to the first hitting times of zero for $X^{(1)}$ and $X^{(2)}$, respectively.
Since $X^{(1)}$ and $X^{(2)}$ are independent, given $X_0=x$, it follows that
$L^{(1)}$ and $L^{(2)}$ are also independent, given $X_0=x$.

Recall (see
 \cite{shigawatanabe73} and  \cite{borodinsalminen02}
p. 72) that
if $Y^{(1)}$ and $Y^{(2)}$ are two independent non-negative
time-homogeneous
diffusions then $Y^{(1)} + Y^{(2)}$
is a diffusion if and only if the generators of $Y^{(1)}$ and $Y^{(2)}$ are of the forms
$$a x \frac{d^2}{d x^2} + (b x + c_i) \frac{d}{d x}, \quad i=1,2,$$
respectively,
where $a>0$, $c_i \geq 0$, $b \in \RR$.
The generator of $Y^{(1)} + Y^{(2)}$ is
$$a x \frac{d^2}{d x^2} + (b x + (c_1+c_2)) \frac{d}{d x}.$$
From Theorem $\ref{theorem4}$ we obtain now the Ray--Knight theorem for $L^{(e)}$.
\begin{thm} \label{ell}
  Conditionally on $X_0=x$,
  $$\{L^{(e)}(y): 0 \leq y \leq x\} \rr \{Z_y^{(4, 2 \mu)}: 0 \leq y \leq x\},$$
  $$\{L^{(e)}(x+y):y \geq 0\} \rr \{Z_y^{(0, 2 \mu)}: y \geq 0\},$$
  where
the process
$Z^{(4, 2 \mu)}$ is started from 0 and the process
 $Z^{(0, 2 \mu)}$ is started from the position of
$Z^{(4, 2 \mu)}$  at time $x$ but otherwise  $Z^{(4, 2 \mu)}$ and $Z^{(0, 2 \mu)}$
are independent.
\end{thm}

\begin{lem}\label{lemtau}
  Let $H_0(L):= \inf \{ y \geq 0: L^{(e)} (X_0 + y) = 0\}$. Then
  the random variables $H_0(L)$ and $L^{(e)}(X_0)$ are exponentially distributed
  with parameters $2\mu$ and $\mu$, respectively.
\end{lem}
\textit{Proof:}
Define
$$M^{(1)}:=\sup\{X_t: \, 0 \leq t \leq d_0\} \;\text{and} \; M^{(2)}:=\sup\{X_t: \, g_0 \leq t \leq 0\}.$$
The distribution of $M^{(i)}$  ($i=1,2$) is given by
(see \cite{borodinsalminen02} p. 14)
$$\pr_x (M^{(i)}<y)=\pr_x(H_0(X) < H_y(X))=
\frac{S(y) - S(x)}{S(y) - S(0)}.$$
Let
$\widetilde M^{(i)}:= M ^{(i)}-x$.
Then
$$\pr_x(\widetilde M^{(i)} <y)= \pr_x(M^{(i)}<x+y)=
\frac{e^{2 \mu (x+y)} -e^{2 \mu x}}{e^{2 \mu (x+y)} - 1},$$
where we used that for $RBM$ with drift
$$S(x)=\frac{1}{2 \mu} (1-e^{2 \mu x}).$$
Hence, by integration,
\begin{eqnarray*}
 \pr(H_0(L) <y) & = & \int_0^\infty 2 \mu e^{-2 \mu x}\, \pr_x (\widetilde M^{(1)} < y)
   \pr_x (\widetilde M^{(2)} < y) \,dx \\
 & = & \int_0^\infty 2 \mu e^{-2 \mu x}
\left( \frac{e^{2 \mu (x+y)} -e^{2 \mu x}}{e^{2 \mu (x+y)} - 1} \right)^2  dx \\
& = & 1-e^{-2 \mu y}.
\end{eqnarray*}

Next we show that $L^{(e)}(X_0) \sim Exp(\mu)$, that is $L^{(e)}(X_0)$ is
exponentially distributed with parameter $\mu$.
Recall that for Brownian motion with drift $- \mu$ killed at zero
$$\widehat G_0 (x,x) = \frac{1}{2 \mu} \left( e^{2 \mu x} -1 \right).$$
The Laplace transform of $L(H_0,x)$  is
(see \cite{borodinsalminen02} p. 32)
\begin{equation} \label{ell1}
\ex_x \left(e^{-\alpha L(H_0(X),\,x)} \right) = \frac{1}{\widehat G_0(x,x)
m(x) \alpha +1} = \frac{\mu}{(1-e^{-2 \mu x}) \alpha + \mu}.
\end{equation}
From $(\ref{ell1})$ we find
$$  \ex \big( e^{-\alpha L^{(e)}(X_0)}\big) =
  \int_0^\infty 2 \mu e^{-2 \mu x} \Big(\frac{\mu}{(1-e^{-2 \mu x})
  \alpha + \mu} \Big)^2 \, dx
  = \frac{\mu}{\mu + \alpha}. $$
  \qed

\begin{prop} \label{rbm}
 The random variables $I^{(+)}$, $I^{(-)}$, and $d_0$ have the same law.
\end{prop}
\noindent
\textit{Proof:} We show first that $I^{(+)} \rr I^{(-)}$. By Proposition $\ref{ell}$ and
$(\ref{occup})$,
\begin{eqnarray*}
  \ex \big(e^{-\alpha I^{(+)}} \big) & = & \ex \Big( \exp \Big(
  -\alpha \int_0^\infty \indic_{ \{X_s > X_0 \}}ds
  \Big) \Big) \\
& = & \ex \Big( \exp \Big( - \alpha \int_0^\zeta Z^{(0, 2 \mu)}_y \, dy
   \Big) \Big),
\end{eqnarray*}
where $Z^{(0, 2 \mu)}_0 \sim Exp(\mu)$  and $\zeta \sim Exp(2 \mu)$ is the
 life time of
$Z^{(0, 2 \mu)}$.
Note that (cf. \cite{borodinsalminen02}
p. 90)
$$\{Z^{(0, 2 \mu)}(\zeta-y): 0 \leq y \leq \zeta \} \rr \{Z^{(4, 2 \mu)} (y):
0 \leq y \leq \tau\},$$
where $\tau \sim Exp(2 \mu)$ is independent of $Z^{(4, 2 \mu)}$.
Consequently, by Proposition
$\ref{ell}$,
\begin{eqnarray*}
  \ex \big(e^{-\alpha I^{(+)}} \big) & = &
  \ex \Big( \exp \Big(-\alpha \int_0^\zeta Z^{(0, 2 \mu)}_y \, dy
   \Big) \Big)\\
   & = & \ex \Big(\exp \Big(-\alpha \int_0^\zeta Z^{(0, 2 \mu)}_{(\zeta -y)} \, dy
   \Big) \Big)\\
   & = & \ex \Big(\exp \Big(-\alpha \int_0^\tau Z^{(4, 2 \mu)}_y \, dy
   \Big) \Big)\\
   & = & \ex \big(e^{-\alpha I^{(-)}}\big).
\end{eqnarray*}
To deduce that $I^{(+)} \rr d_0$, we adopt the argument in 
\cite{salminenyor04} used in the proof of the Biane--Imhof identity.
First recall that
$$I^{(+)} = \int_0^\infty Z_s \, d s,$$
where the process $Z= Z^{(0,2 \mu)}$ is the solution of the SDE
\begin{equation} \label{sde}
  dZ_t = 2 \sqrt{Z_t}\, d W_t -2 \mu\, Z_t \, d t, \quad Z_0 \sim Exp(\mu).
\end{equation}
Introduce next
$$A_t:= \int_0^t Z_s \, ds$$
and let $\alpha_t$ be the right-continuous inverse of $A_t$.
Since the quadratic variation of the local martingale $Y_t = \int_0^t \sqrt{Z_s} \, dW_s$ is given by
$$\langle Y, Y\rangle_t = \int_0^t Z_s \, ds = A_t,$$
it follows from Lévy's characterization theorem that $Y_{\alpha_t}$ is a Brownian
motion, say $B_t$, started at zero and stopped at $A_\infty$. Hence, for $t < A_\infty$,
\begin{eqnarray*}
  Z_{\alpha_t} - Z_0 & = & 2 \int_0^{\alpha_t} \sqrt{Z_s} \, d W_s -
  2 \mu \int_0^{\alpha_t} Z_s \, ds \\
  & = & 2 Y_{\alpha_t} - 2 \mu A_{\alpha_t} \\
  & = & 2 B_t - 2 \mu t.
\end{eqnarray*}
Letting $t \to A_\infty$ gives
$Z_0 / 2 = B^{(\mu)}_{A_\infty},$
where $B^{(\mu)}_t := -B_t + \mu t$. Thus, taking into account that
$$0< Z_{\alpha_t} = Z_0+2 B_t -2 \mu t, \quad 0 \leq t < A_\infty,$$
yields
\begin{equation} \label{ainf}
A_\infty = \inf\{t: B^{(\mu)}_t = Z_0/2\}.
\end{equation}
Since $Z_0 \sim Exp(\mu)$  and is
independent of $B^{(\mu)}$, $(\ref{ainf})$ is equivalent with
\begin{equation} \label{ainf1}
I^+ = A_\infty \rr H_\lambda (B^{(\mu)}),
\end{equation}
where $\lambda \sim Exp(2 \mu)$ is independent of $B^{(\mu)}$. Noting that the right-hand
side of $(\ref{ainf1})$ is identical in law to $H_0 (X)$ gives $I^{(+)} \rr d_0$,
as claimed.
\qed

\section{Relation to It\^{o} excursion theory} \label{excursions}

\subsection{Distribution of $(-g_0, d_0)$ in terms of Lévy measure} \label{subinverse}

Consider the distribution of $d_0$ (and $-g_0$). From $(\ref{jointlap})$,
\begin{equation} \label{lapd}
\ex\left( e^{-\alpha d_0} \right) = \ex\left( e^{\alpha g_0} \right)
= \frac{1}{M \alpha G_\alpha(0,0)}.
\end{equation}
Next let $\ell(t,0)$ be the local time at zero up to time $t$ (with respect to the speed
measure) of $\{X_s : s \geq 0\}$, $X_0=0$.
Let $A= \{A_t : t \geq 0\}$ be the right-continuous inverse of $\ell$.
 The process $A$ is
a subordinator and
$$\ex_0 \big(\exp (-\alpha A_t) \big) = \exp (- t\, \Psi(\alpha)),$$
where the Laplace exponent $\Psi$
is given by
$$\Psi(\alpha) = \int_0^\infty \left(1-
e^{-\alpha t}\right) n^+(dt) =\alpha \int_0^\infty e^{-\alpha t} n^+(t, \infty)\, dt$$
with the associated Lévy measure $n^+(dt)$.

\begin{prop} \label{lapexp}
  The Laplace exponent $\Psi(\alpha)$ of $A$ is given by
  \begin{equation} \label{green00}
  \Psi(\alpha) = \frac{1}{G_\alpha(0,0)}
  = - \frac{d}{d S(x)} \ex_x \left(e^{-\alpha H_0} \right) \Big{|}_{x=0+}.
\end{equation}
\end{prop}
\noindent
\textit{Proof:} See  \cite{itomckean74} p. 214. \qed

\begin{prop} \label{propdglevy}
  Let $n^+(dt)$ be the Lévy measure of $A$. Then
\begin{equation} \label{l1}
    \pr(-g_0 \in dt)= \pr(d_0 \in dt) = \frac{n^+(t,\infty)}{M}\, dt;
\end{equation}
\begin{equation} \label{l2}
 \pr(V \in dv) = \frac{v}{M} \, n^+(dv), \; \text{where} \; V:= d_0-g_0;
\end{equation}
\begin{equation} \label{dglevy}
  \pr(d_0>v, -g_0>w)=\frac{1}{M} \int_{v+w}^\infty n^+(t, \infty)\, dt;
\end{equation}
 \begin{equation} \label{joint1}
\pr(d_0 \in dt, -g_0 \in ds) / dt \, ds = -\frac{1}{M} \frac{d}{d t} n^+(t+s, \infty).
\end{equation}
\end{prop}
\noindent
\textit{Proof:} Formula $(\ref{l1})$
is a consequence of $(\ref{lapd})$ and Proposition $\ref{lapexp}$, and
is given in \cite{pitman86} and  \cite{pitmanyor97, pitmanyor03}
(the \it global formula\rm).
Formulae $(\ref{l2})$, $(\ref{dglevy})$, and $(\ref{joint1})$ follow immediately from
$(\ref{l1})$ and Proposition $\ref{propk}$. \qed

\begin{rema} \label{remal} \rm
It is interesting to compare the right-hand sides of $(\ref{joint1})$ and $(\ref{joint})$.
For this, notice first from $(\ref{green00})$ that
\begin{equation} \label{ninfty}
    n^+(t, \infty)= \frac{d}{dS(x)} \pr_x (H_0 \geq t) \Big{|}_{x=0+}.
\end{equation}
Thus, we have
\begin{eqnarray*}
  -\frac{d}{dt} n^+(t+s, \infty) & = & \frac{d}{dt} \frac{d}{dS(x)}
    \pr_x (H_0<t+s) \Big{|}_{x=0+} \\
& = & \frac{d}{dS(x)}\frac{d}{dt}\pr_x (H_0<t+s) \Big{|}_{x=0+} \\
& = & \frac{d}{dS(x)} \frac{d}{dS(y)} \widehat p(t+s;x,y) \Big{|}_{x,y=0+}.
\end{eqnarray*}
\end{rema}

\subsection{Spectral representations for $d_0$ and $V=d_0-g_0$} \label{subspectral}

In this section we show that the distribution of $d_0$ (and $-g_0$) is a mixture of
exponential distributions and the distribution of $V = d_0-g_0$ is a mixture of
gamma distributions. The mixing measures are the same and closely related to the
so called principal spectral measure of the process $X$, as defined in
Krein's
theory of strings, see  \cite{kackrein74, kotaniwatanabe81, kuchler86}.
Our starting point is the result due to Knight (see \cite{knight81}) which states that
the Lévy measure $n^+(dt)$ is absolutely continuous with respect to the Lebesgue
 measure and there exists
 a unique measure $\Delta$ such that
\begin{equation} \label{krein}
\nu(t):= n^+(dt)/ dt  = \int_0^\infty e^{-z t} \,\Delta(dz).
\end{equation}
Moreover, $\Delta$ has the properties
\begin{equation} \label{rep}
   \int_0^\infty \frac{\Delta(dz)}{z(z+1)} < \infty
\end{equation}
and
\begin{equation} \label{repdelta}
   \int_0^\infty \frac{\Delta(dz)}{z} = \infty.
\end{equation}
We remark (cf. \cite{knight81}) that $(\ref{rep})$ is equivalent with the
defining property of the Lévy measure of a subordinator, i.e.,
$$\int_0^\infty (1 \wedge t) \, n^+(dt) < \infty.$$
The property in $(\ref{repdelta})$ is a consequence of the fact that the speed measure
is strictly positive everywhere and, in particular,
in a right neighbourhood  of $0$, see \cite{kackrein74} p. 82.
For another proof  of $(\ref{repdelta})$ see
 \cite{kuchlersalminen89}, where $\Delta$ is interpreted
as the principal spectral measure of a killed string.
\begin{prop} \label{propkrein}
 Let $\Delta$ be the measure introduced above
  and associated with the process $X$.Then the measure
 $$\widetilde \Delta(dz)=\Delta(dz) / (M z^2)$$
 is a probability measure. Moreover,
  \begin{equation} \label{rep1}
     \pr(d_0 \in dt) / dt = \int_0^\infty z e^{-z t}\, \widetilde\Delta(dz),
  \end{equation}
 and
  \begin{equation} \label{rep2}
     \pr(V \in dv) / dv = \int_0^\infty z^2 v \,e^{-zv}\, \widetilde \Delta(dz).
  \end{equation}
\end{prop}
\noindent
\textit{Proof:}
Recall that in the recurrent case (see \cite{salminen93} p. 220 and
\cite{borodinsalminen02} p. 20),
\begin{equation}\label{rec}
 \lim_{\alpha \searrow 0} \alpha G_\alpha(x,x) = \frac{1}{m\{I\}}, \quad \text{for
  all } x \in I.
\end{equation}
By $(\ref{rec})$, $(\ref{krein})$, and Fubini's theorem,
\begin{eqnarray*}
 M:=m\{I\} &=& \lim_{\alpha \searrow 0} \frac{1}{\alpha G_\alpha(0,0)} \\
 &=& \int_0^\infty  n^+(t, \infty) \, dt \\
 &=&\int_0^\infty dt \int_t^\infty \nu(s) \, ds\\
 &=& \int_0^\infty dt \int_0^\infty \Delta(dz)
 \frac{e^{-zt}}{z} \\
 &=& \int_0^\infty \frac{\Delta(dz)}{z^2},
\end{eqnarray*}
and, therefore, $\widetilde \Delta$ is a probability measure.
Formulae $(\ref{rep1})$ and $(\ref{rep2})$ follow now from Proposition $\ref{propdglevy}$ and spectral representation
   $(\ref{krein})$.  \qed

\begin{rema}  \label{recur} \rm
From the proof of Proposition $\ref{propkrein}$ a new test for positive
 recurrence emerges:
   a recurrent diffusion $X$ is positively recurrent if and only if
   $$\int_0^\infty \frac{\Delta(dz)}{z^2} < \infty.$$
\end{rema}

\subsection{Excursion bridges} \label{palm}

In this section we use the theory of the Palm measures
to show that the result in Theorem \ref{theorem1} has a
counterpart in the framework of excursions from 0 to 0 of
the stationary non-negative diffusion $X.$
Let
$$
\mathcal{M}:= \{t \in \RR: X_t = 0\}
$$
be the zero set of $X$
and
$$
\mathcal{L}:= \{t\in{\cal M}\,:\, \exists\ \varepsilon>0\ \forall\
0<s<\varepsilon\quad   X_{t+s} > 0\}
$$
be the set of the starting times of excursions from 0 to 0 on a path of
$X.$

We introduce next the space $E$ consisting of continuous
functions $e:\RR_+\mapsto\RR_+$ such that $e(0)=0$ and for which there
exists $ \zeta=\zeta(e)>0$ with the property $e(t)>0$ for
$0<t<\zeta$ and $e(t)=0$  for $t\geq \zeta.$
The space $(E,\cal E),$ where $\cal E$ is the $\sigma$-algebra
generated by the cylinder sets in the usual way,
is called the canonical excursion
space for excursions from 0 to 0 (associated with $X$).
For $t\in{\cal L}$ define $X^{(ex,t)}=\{X^{(ex,t)}_{s}\,:\, s\geq 0\}$ where
$$
X^{(ex,t)}_{s}:=
\begin{cases}
X_{t+s},& $for\ $ t+s<R\\
0,& $for\ $ t+s\geq R,\\
\end{cases}
$$
and $R:=\inf\{u>t\,:\, X_u=0\}.$ Clearly, $X^{(ex,t)}\in E$ and
$\zeta(X^{(ex,t)})=R.$

We follow now  \cite{pitman86}
and define the concept of
equilibrium excursion measure which is a particular case of the Palm
measure construction. See \cite{pitman86}
for general results on Palm measures and further references.
For background on Palm measures, especially in queueing,  see,
e.g., 
\cite{baccellibremaud03}.

\begin{defi}
\rm For $B\in{\cal E}$ let
\begin{equation} \label{expalm}
\prq(B):= \ex\left( \left|\{t: 0 <t <1, t \in {\cal L}, X^{(ex,t)} \in B\}\right|\right),
\end{equation}
where $|\cdot|$ denotes the number of elements in the set
under the consideration. The measure $\prq$ is called the equilibrium
excursion measure for the excursions of $X$ from 0 to 0.
\end{defi}

The next proposition and its corollary are adopted from  \cite{pitman86} (Theorem p. 290) where a more general
statement originating from  \cite{mecke67} and \cite{neveu77} is given. Because of the key
importance of this result for our
application, we restate and formulate it especially for excursions
(cf. Section III in \cite{pitman86} where
connections with the Maisonneuve formula are discussed).

\begin{prop}
\label{palmproperties}
The measure $\prq$ is $\sigma$-finite and for
all measurable $f: \RR \times E \to [0, \infty)$,
 \begin{equation}
 \label{palm1}
    \ex \Big( \sum_{t \in {\cal L}}f(t, X^{(ex,t)}) \Big) = \int_\RR \int_E ds \,
   \prq(d e)\,f(s, e).
\end{equation}
\end{prop}
\textit{Proof:}
For the proof of the first claim, see  \cite{pitman86}.
To prove formula $(\ref{palm1})$, observe from
the definition of $\prq$, that $(\ref{palm1})$ is equivalent with
\begin{equation} \label{p1}
 \ex \Big( \sum_{t \in {\cal L} }f(t, X^{(ex,t)}) \Big) =
 \ex \Big( \sum_{t \in {\cal L}\cap (0,1)} \int_\RR ds f(s, X^{(ex,t)})\Big).
\end{equation}
Let $\theta=\{ \theta_t: t \in \RR\}$ be the usual shift operator defined in the underlying
probability space
via
$$
X_s\circ\theta_t(\omega) =X_{t+s}(\omega)\quad\forall t,s\in\RR.
$$
Notice that the law of $X$ is invariant under $\theta$; this is, in our case, equivalent with
 the stationarity of $X$. Consequently, by stationarity,
\begin{eqnarray*}
&&\ex \Big( \sum_{t \in {\cal L}\cap (0,1)}
\int_\RR ds\, f(s, X^{(ex,t)})\Big)\\
&&\hskip2cm
=
\int_\RR ds\
 \ex \Big(   \sum_{t \in {\cal L} \circ \theta_s\cap(0,1)} f(s, X^{(ex,t)}\circ\theta_s) \Big)\\
&&\hskip2cm
=
\int_\RR ds\
 \ex \Big( \sum_{t \in (0,1), t+s \in {\cal L}}
   f(s, X^{(ex,t+s)})\Big)\\
&&\hskip2cm
= \int_0^1 dt\
 \ex \Big( \sum_{s \in \RR, s \in {\cal L}\circ\theta_t}
   f(s, X^{(ex,s)}\circ\theta_t)\Big) \\
&&\hskip2cm
= \int_0^1 dt\
 \ex \Big(  \sum_{s \in {\cal L}}
   f(s, X^{(ex,s)})\Big) \\
&&\hskip2cm
=  \ex \Big( \sum_{s \in {\cal L}}
   f(s, X^{(ex,s)})\Big), \\
\end{eqnarray*}
hence $(\ref{p1})$ is proved giving $(\ref{palm1})$.
\qed

\begin{cor}
\label{cor}
For a jointly measurable $h:\RR\times E\mapsto [0,\infty)$
\begin{equation} \label{iii2}
\ex\left( h(-g_0, X^{(ex,g_0)}) \right)
 = \int_E\prq(de)\Big(\int_0^{\zeta(e)}
 h(s,e) \, ds\Big).
\end{equation}
In particular, for  $a>0,$ $v>0,$ and $e\in E$
\begin{equation} \label{iii1}
  \pr(-g_0 \in da, \, X^{(ex,g_0)}\in de) = da \, \prq(de) \,
  \indic_{(a<\zeta(e))},
\end{equation}
\begin{equation} \label{iii111}
  \pr(X^{(ex,g_0)}\in de) = \prq(de) \zeta(e),
\end{equation}
\begin{equation} \label{iii112}
    \pr(-g_0 \in da)/da = \prq(\zeta>a),
\end{equation}
and
\begin{equation} \label{iii11}
  \pr(V \in dv)= v\,\prq(\zeta(e)\in dv).
\end{equation}
\end{cor}
\textit{Proof:}
To prove formula (\ref{iii2}), we set in (\ref{palm1})
$$
f(t, e) = h(-t, e) \indic_{(0<-t<\zeta(e))}
$$
with $h$ jointly
measurable (cf. \cite{pitman86} p. 291 and \cite{neveu77} p. 332). Then,
due to the particular form of $f,$
\begin{eqnarray}
\label{iii}
\nonumber
&&\ex \Big(  \sum_{t \in {\cal L}}
   f(t, X^{(ex,t)})\Big)
 =\ex\ \left( h(-g_0, X^{(ex,g_0)}) \right)\\
&&\hskip3.3cm
 = \int_E\prq(de)\Big(\int_0^{\zeta(e)}
 h(s,e) \, ds \Big).
\end{eqnarray}
Choosing $h$ appropriately in (\ref{iii2}) yields (\ref{iii1}),
(\ref{iii111}), and (\ref{iii112}). Finally,  take
in (\ref{iii}) $h(s,e)={\bf 1}_B(e),\ B=\{\zeta>v\},\ v>0,$
to obtain
$$
\pr\left( V>v \right)
 = \prq\left(\zeta(e)\,;\,\zeta(e)>v\right)
$$
which is equivalent with (\ref{iii11}).
\qed

\begin{rema} \label{palmlevy} \rm
1)\,From (\ref{iii1}) and  (\ref{iii111}) it follows that
the distribution of $-g_0$ given the excursion is uniform on $(0,V)$ and, in particular,
 $(-g_0, d_0) \in \kk$
(cf. Proposition $\ref{proposition1}$).\\
2)\,It is seen from Proposition $\ref{propdglevy}$ and Corollary $\ref{cor}$ that
 $$\prq (\zeta(e) \in dv) = \frac{n^+ (dv)}{M}.$$
\end{rema}



To proceed, let ${\bf M}$ denote the It\^{o} excursion law for the
excursions of $X$ from 0 to 0. Then it is well-known that
$$
{\bf M}(\zeta>a)=c\,\int_0^\infty m(dx)\,n_x(0,a),
$$
where $c$ is a normalizing constant and
$$
n_x(0,a):=\pr(H_0\in da)/da
$$
is the density of the first hitting time $H_0:=\inf\{t\,:\, X_t=0\}$
given that $X_0=x.$ Now (\ref{iii112}) yields
$$
{\bf M}(\zeta>a)=c' \prq(\zeta>a)
$$
for some (normalizing) constant $c'.$ In fact, as can be deduced from
(\ref{iii1}), the measures $\bf M$ and $\prq$ are the same (up to a
constant):
\begin{prop}
\label{mq}
There exists a constant $c'$ such that for all $B\in{\cal E}$
\begin{equation}
\label{ito}
{\bf M}(B)=c' \prq(B).
\end{equation}
\end{prop}
We refer to \cite{pitman86} for a discussion about (\ref{ito})
and the associated normalizations. See also  \cite{bismut85}
for the Brownian motion case.

To formulate the main result of this section, we need the concept of
excursion bridge of $X.$ For reflecting Brownian motion the excursion
bridge is a 3-dimensional Bessel bridge from 0 to 0 (of some length
$l$). In general, the excursion bridge of $X$ from 0 to 0 of the
length $l$  can be described as the
process $X^{(0,l,0)}$ obtained from  $X$ with $X_0=x>0$ conditioned
to hit 0 for the first
time at time $l$ by letting $x\to 0.$ Clearly, defining
$X^{(0,l,0)}_t=0$ for $t\geq l$ it holds $X^{(0,l,0)}\in E.$

\begin{thm} \label{propex}
Let $X^{(0,l,0)}=\{X^{(0,l,0)}_t\,:\,t\geq 0\}$ be the excursion
bridge as defined above and let $U$ be a uniformly on $(0,l)$  distributed
random variable independent of $X^{(0,l,0)}.$ Let
$$
I^{(l,+)}:=\int_0^l{\bf 1}_{\{X^{(0,l,0)}_s>X^{(0,l,0)}_U \}}\,ds
$$
and
$$
I^{(l,-)}:=\int_0^l{\bf 1}_{\{X^{(0,l,0)}_s<X^{(0,l,0)}_U\}}\,ds.
$$
Then $I^{(l,+)}$ and $I^{(l,-)}$ are identically and uniformly on $(0,l)$
distributed random variables.
\end{thm}
\textit{Proof:}
We use formula (\ref{iii2}) and let therein
$$
h(s,e)=F(\zeta(e))\,G \Big(\int_0^{\zeta(e)}{\bf
  1}_{\{e_t>e_s\}}dt \Big)\,{\bf 1}_{\{s<\zeta(e)\}}
$$
with  Borel-measurable functions  $F$ and $G.$
Recall the notation $V=d_0-g_0,$ and
consider first the left-hand side of
(\ref{iii2}) with $h$ as above:
\begin{eqnarray}\label{left}
\nonumber
&&\ex\left(h(-g_0,X^{(ex,g_0)})\right)
=\ex\Big(F(V)\,G \Big(\int_0^{V}{\bf
  1}_{\{X^{(ex,g_0)}_t>X^{(ex,g_0)}_{-g_0}\}}dt \Big)\,{\bf
  1}_{\{-g_0<V\}}\Big)\\
\nonumber
&&\hskip3.5cm
=\ex\Big(F(V)G \Big(\int_0^{V}{\bf
  1}_{\{X_{g_0+t}>X_0\}}dt \Big)\Big)
\\
&&\hskip3.5cm
=\ex\left(F(V)\,G(I^{(+)})\right),
\end{eqnarray}
where, as before,
$$
I^{(+)}:=\int_{g_0}^{d_0}{\bf
  1}_{\{X_{t}>X_0\}}dt.
$$
For the right-hand side of (\ref{iii2}) we have
\begin{eqnarray}
\label{right1}
\nonumber
&&
\int_E\prq(de)\Big(\int_0^{\zeta(e)}
 h(s,e) \, ds\Big)\\
\nonumber
&&\hskip0.5cm
 =
\prq\Big(F(\zeta(e))\,\zeta(e) \int_0^{\zeta(e)}
G \Big(\int_0^{\zeta(e)}{\bf
  1}_{\{e_t>e_s\}}dt \Big)\, \frac{ds}{\zeta(e)}\Big)
\\
&&\hskip0.5cm
=
\frac{1}{c'}\,{\bf M}\Big(F(\zeta(e))\,\zeta(e) \int_0^{\zeta(e)}
G \Big(\int_0^{\zeta(e)}{\bf
  1}_{\{e_t>e_s\}}dt \Big)\, \frac{ds}{\zeta(e)}\Big),
 \end{eqnarray}
where in the last step Proposition \ref{mq} is used. Next recall that the description of the It\^{o}
measure via the lengths of the excursions (see, e.g., 
\cite{revuzyor01} p. 497,  \cite{borodinsalminen02} p. 60, \cite{rogerswilliams87} p. 421) says that
 for all $B\in{\cal E}$,
$$
{\bf M}(B)=\int_0^\infty\,{\bf M}(\zeta\in dl)\, \pr^{(0,l,0)}(B),
$$
where $\pr^{(0,l,0)}$ is the probability measure associated with the excursion bridge process
$X^{(0,l,0)}$ defined in the canonical excursion space $E.$ Consequently, we obtain
 \begin{eqnarray}
\label{right2}
\nonumber
&&{\bf M}\Big(F(\zeta(e))\,\zeta(e) \int_0^{\zeta(e)}
G\Big(\int_0^{\zeta(e)}{\bf
  1}_{\{e_t>e_s\}}dt\Big)\, \frac{ds}{\zeta(e)}\Big)
\\
\nonumber
&&\hskip.5cm
=
\int_0^\infty\,{\bf M}(\zeta(e)\in dl)\, F(l)\, l\ \ex^{(0,l,0)}\Big(
\int_0^{l}
G\Big(\int_0^{l}{\bf
  1}_{\{e_t>e_s\}}dt\Big)\, \frac{ds}{l}\Big)
\\
&&\hskip.5cm
=
\int_0^\infty\,{\bf M}(\zeta(e)\in dl)\, F(l)\, l\ \ex^{(0,l,0)}\Big(
G\Big(\int_0^{l}{\bf
  1}_{\{e_t>e_U\}}dt\Big)\Big),
\end{eqnarray}
where $U$ is uniformly on $(0,l)$ distributed random variable independent
of $X^{(0,l,0)}.$ By (\ref{iii2}), the right-hand sides of (\ref{left})
and (\ref{right1}) are equal. Consequently, using (\ref{iii11}) in Corollary \ref{cor}
and (\ref{right2}), we obtain
$$ \int_0^\infty \pr(V \in dl)  F(l) \, \ex \left( G(I^{(+)}) | V=l \right)
= \int_0^\infty \pr(V \in dl)  F(l) \, \ex^{(0,l,0)} \left( G(I^{(l,+)})\right)$$
giving
$$\ex \left( G(I^{(+)}) | V=l \right) = \ex^{(0,l,0)} \left( G(I^{(l,+)})\right).
$$
Combining this with the result in Theorem \ref{theorem1} completes
the proof.
\qed

\section{Remarks on the null recurrent case}  \label{nullrec}

\begin{rema} \label{nullrec1}
\rm
Suppose now that $X$ is null recurrent.
Then  the speed measure $m(dx)$ serves still as the stationary measure
of $X$ but because $m\{I\} =\infty$ it cannot be normalized to a probability measure.
However, the result in Theorem $\ref{theorem1}$ is also valid in this case:
\begin{eqnarray} \label{eqnull}
 \int_I m(dx) \, \ex_x \big( e^{-\alpha I^{(+)} - \beta I^{(-)}}\big)
 &=& \int_I m(dx) \, \ex_x \big(e^{-\alpha d_0 + \beta g_0} \big)  \nonumber \\
 &=& \frac{1}{\alpha-\beta} \left(  \frac{1}{G_\alpha(0,0)} -\frac{1}{G_\beta(0,0)} \right).
\end{eqnarray}
In particular, given that $V:=d_0-g_0=v$, $I^{(+)}$ and $I^{(-)}$ are uniformly
distributed on $(0,v)$.
In the case when $X$ is a stationary reflecting Brownian motion living above zero
$(\ref{eqnull})$ gives
\begin{eqnarray*}
 \int_I m(dx) \, \ex_x \big( e^{-\alpha I^{(+)} - \beta I^{(-)}}\big)
 &=& \int_I m(dx) \, \ex_x \big(e^{-\alpha d_0 + \beta g_0} \big) \\
 &=& \frac{1}{\alpha-\beta} \left(  \sqrt{2 \alpha} -\sqrt{2 \beta} \right) \\
 &=& \frac{\sqrt{2}}{ \sqrt{\alpha} + \sqrt{\beta}}.
\end{eqnarray*}

It is also easily checked that finiteness of $m$ is not needed
in the proof of Theorem $\ref{propex}$ and, hence, for all diffusion excursion bridges
$X^{(0,l,0)}$ of length $l$ the random variables $I^{(l,+)}$ and $I^{(l,-)}$
(defined as in Theorem $\ref{propex}$) are uniformly distributed on $(0,l)$.
\end{rema}

\begin{rema} \label{vervaat}
\rm
For a Brownian excursion there is an alternative simple proof of the result
in Theorem $\ref{propex}$ based on Vervaat's path transformation
(see \cite{vervaat79, biane86, bertoin03}).
Indeed, let $X^{(0,l,0)}$ denote a standard Brownian excursion of length $l$, that is
a 3-dimensional Bessel bridge of length $l$, and $U$ be
a random variable uniformly distributed on $(0,l)$ and independent of $X^{(0,l,0)}$.
Then the process $X^{br} = \{ X^{br}_t: 0 \leq t \leq l\}$ defined as
$$X^{br}_t = \left \{
    \begin{array}{ll}
       X^{(0,l,0)}_{t+U} - X^{(0,l,0)}_{U}, & t + U \leq l, \vspace{0.3cm}\\
       X^{(0,l,0)}_{t+U-l} - X^{(0,l,0)}_{U}, & t + U \geq l
     \end{array}
     \right.$$
  equals in distribution to a standard Brownian bridge of length $l$.
Clearly, with the notation in Theorem $\ref{propex}$,
$$I^{(l,+)} = \int_0^l \indic_{\{X^{br}_s >0\}} \, ds =: I^{(br,+)}$$
and
$$I^{(l,-)} = \int_0^l \indic_{\{X^{br}_s <0 \}} \, ds =: I^{(br,-)}.$$
The claim follows now from the well-known
result due to Lévy, saying that $I^{(br,+)}$ and $I^{(br,-)}$ are
uniformly distributed
on $(0,l)$ (see, e.g., \cite{borodinsalminen02} p. 163 and  \cite{yor95}
p. 43).
\end{rema}

\section{Examples} \label{subexamples}

To illustrate the results in Section $\ref{results}$ we consider some examples of
stationary diffusions
in stationary state.

\begin{exa} \label{example1} \bf (Reflecting Brownian motion with drift)
\rm Let $X$ be a stationary reflecting Brownian motion
with drift $-\mu<0$ living above zero. The speed measure is
$$m(dx)= 2 e^{-2 \mu x} \, dx,$$
the scale function is
$$S(x)=\frac{1}{2 \mu}(1- e^{2 \mu x}),$$
and the Green function at $(0,0)$ is
$$G_\alpha(0,0)=\frac{1}{\sqrt{2 \alpha + \mu^2} -\mu}.$$
Hence by Theorem $\ref{theorem1}$ and $(\ref{jointlap})$
(cf. \cite{salminennorros01}),
\begin{eqnarray}
\label{laprefl}
  & & \ex \big( e^{-\alpha I^{(+)} -\beta I^{(-)}} \big) =
  \ex \left( e^{-\alpha d_0 + \beta g_0} \right) \nonumber \\
  & & \hspace{1.2cm} = \frac{\mu}{\alpha - \beta} \left( (\sqrt{2 \alpha + \mu^2}-\mu) -
  (\sqrt{2 \beta + \mu^2}-\mu)\right) \nonumber \\
  & & \hspace{1.2cm} = \frac{2 \mu}{\sqrt{2 \alpha +\mu^2} + \sqrt{2 \beta + \mu^2}}.
\end{eqnarray}
The transition density for the killed process is
$$\widehat p(t; x,y) = \frac{1}{2 \sqrt{2 \pi t} } e^{\mu(x+y)-\frac{\mu^2 t}{2}}
\left(e^{-\frac{(x-y)^2}{2t}} -e^{-\frac{(x+y)^2}{2t}}\right).$$
Thus by $(\ref{joint})$, the joint density of  $(-g_0,d_0)$ and $(I^{(+)}, I^{(-)})$ is
\begin{eqnarray*}
f(t,s) & = & \mu
\left ( \frac{d}{d x} \frac{d}{d y}
    \widehat p(s+t; x, y) /
    (S'(x) S'(y)) \right) \Big|_{x, y=0}   \\
    & = & \frac{\mu}{\sqrt{2 \pi (t+s)^3}} e^{-\frac{\mu^2}{2} (t+s)}.
\end{eqnarray*}
To check the obtained formulae  we give here also the following
expressions for the corresponding Laplace transforms when $X$ started at $x$:
\begin{eqnarray*}
& & \ex_x \left( \exp \Big( -\alpha \int_0^{H_0(X)} \indic_{ \{ 0 \leq X_s < x \}} \, ds
- \beta \int_0^{H_0(X)} \indic_{\{X_s > x\}} \, ds\Big)\right) \\
& & \hspace{0.5cm} \frac{\sqrt{2 \alpha + \mu^2} e^{\mu x}}{ \sqrt{2 \alpha +
\mu^2} \cosh(x \sqrt{2 \alpha + \mu^2}) + \sqrt{2 \beta + \mu^2} \sinh(
x \sqrt{2 \alpha + \mu^2})}
\end{eqnarray*}
(see \cite{borodinsalminen02}, formula 2.2.6.1 p. 300),
while
$$\ex_x \left( e^{-\alpha H_0(X) }\right) = e^{(\mu-\sqrt{2 \alpha + \mu^2})x}.$$
\end{exa}

\begin{exa} \label{exa2} \bf (Brownian motion reflected at $0$ and $1$)
 \rm Consider now  a stationary Brownian motion $X$ living in the interval $I=[0,1]$
 and reflected at $0$ and
 at $1$. The speed measure of $X$ is $m(dx)=2 \,dx$, the scale function
 $S(x) = x$, and the Green function at $(0,0)$ is
(see \cite{borodinsalminen02}, p. 122)
  $$G_\alpha(0,0)=\frac{\coth(\sqrt{2 \alpha})}{\sqrt{2 \alpha}}.$$
   Hence the Laplace transform of $(-g_0, d_0)$ and $(I^{(+)},I^{(-)})$
   is  ($\alpha \ne \beta$)
  \begin{eqnarray*}
   & & \ex \big( e^{-\alpha I^{(+)} - \beta I^{(-)}}\big) =
    \ex \big(e^{-\alpha d_0 +\beta g_0} \big) \\
  & & \hspace{1cm} = \frac{1}{\alpha-\beta} \left( \sqrt{\alpha / 2} \tanh(\sqrt{2 \alpha})
    - \sqrt{\beta / 2} \tanh(\sqrt{2 \beta}) \right)
  \end{eqnarray*}
  and
  $$\ex (e^{-\alpha (d_0-g_0)}) = \frac{1}{4 \sqrt{2 \alpha} \cosh^2(\sqrt{2 \alpha})}
  \left(\sinh(2 \sqrt{2 \alpha})+2 \right). $$

\end{exa}

\begin{exa} \label{exa4} \bf (Squared radial Ornstein-Uhlenbeck process)
 \rm Let $X=\{X_t : t \in \RR\}$ be a stationary squared radial Ornstein-Uhlenbeck process with
 parameters $\nu= n/2-1$ and $\gamma$. The generator of $\{X_t: t \geq 0\}$ is
 $$\mathcal{G} = 2 x \,\frac{d^2}{dx^2} + (n -2 \gamma \,x) \frac{d}{dx}.$$
 The process $X$ can be defined as
 \begin{equation} \label{ou}
 X(t) = e^{-\gamma t} SQBES(e^{2 \gamma t} / 2 \gamma), \quad  t \in \RR,
 \end{equation}
 where $SQBES$ is a squared Bessel process of dimension $n= 2 \nu + 2$ (see
 \cite{pitmanyor97}).
 Assume that $-1 < \nu <0$ and $0$ is a reflecting boundary. Then the process is positively
 recurrent with the speed measure given by (see  \cite{borodinsalminen02} p. 140)
 $$m(dx) = \frac{1}{2} x^\nu e^{-\gamma x} \, dx,$$
 that is (after normalization) the gamma-density with parameters $\gamma$ and $\nu+1$.
 The Green function at $(0,0)$ is given by (see  \cite{borodinsalminen02}
 p. 141)
 $$G_\alpha (0,0) = \gamma^\nu \,\frac{\mathrm{B}(\frac{\alpha}{2 \gamma}, -\nu)}{\Gamma(1+ \nu)},$$
 where $\Gamma(x)$ and $\mathrm{B}(x,y)$ denote the gamma and beta functions, respectively.
 By Theorem  $\ref{theorem1}$,
 \begin{eqnarray*}
   \ex \big(e^{- \alpha I^{(+)} - \beta I^{(-)}}\big) = \ex \left(e^{-\alpha d_0 + \beta g_0}\right)
    & = &\frac{1}{M(\alpha -\beta)}
   \Big( \frac{1}{G_\alpha(0,0)}-\frac{1}{G_\beta(0,0)} \Big)\\
   & = &\frac{2 \gamma}{\alpha -\beta} \Big( \frac{1}{\mathrm{B}(\frac{\alpha}
   {2 \gamma}, - \nu)} -  \frac{1}{\mathrm{B}(\frac{\beta}{2 \gamma}, - \nu)}\Big).
\end{eqnarray*}
Letting $\beta = 0$ and taking the inverse Laplace transform (see  \cite{erdelyi54}
p. 261) gives the density of $d_0$ (and $I^{(\pm)}$) and the expression for $n^+(t, \infty)/ M$:
\begin{equation} \label{example4}
  \pr(d_0 \in dt) / dt = \frac{n^+(t, \infty)}{M} = \frac{2 \gamma}{\Gamma(-\nu) \Gamma(1+ \nu)} e^{2 \mu \nu t}
  \left( 1- e^{-2 \mu t} \right)^\nu.
\end{equation}
 Formula $(\ref{example4})$ can be also obtained using
 $(\ref{ou})$ as in \cite{pitmanyor97}.
The joint density of $d_0$ and $g_0$ (and $I^{(\pm)}$) is given then by
$$f(t,s) = - \frac{1}{M} \frac{d}{dt} n^+(t+s, \infty) = - \frac{4 \mu \nu \gamma}{\Gamma(-\nu)
\Gamma(1+ \nu)} e^{2 \mu \nu (s+t)} \left( 1 - e^{-2 \mu (s+t)}\right)^{\nu -1}.$$

\end{exa}

\bibliography{art1}

\bibliographystyle{plain}

\end{document}